\documentclass[a4paper,12pt,reqno]{amsart}
\usepackage{amsmath, amssymb,amsthm}
\usepackage{enumerate}
\usepackage{cite}

\nonstopmode \numberwithin{equation}{section}
\newtheorem{theorem}{Theorem}[section]

\newtheorem{corollary}{Corollary}[section]
\newtheorem{lemma}{Lemma}[section]
\newtheorem{remark}{Remark}[section]
\newtheorem{problem}{Problem}[section]

\begin{document}
\bibliographystyle{amsplain}

\title{{{
Inclusion properties of Generalized Integral Transform using Duality Techniques
}}}

\author{
Satwanti Devi
}
\address{
Department of  Mathematics  \\
Indian Institute of Technology, Roorkee-247 667,
Uttarkhand,  India
}
\email{ssatwanti@gmail.com}

\author{
A. Swaminathan
}
\address{
Department of  Mathematics  \\
Indian Institute of Technology, Roorkee-247 667,
Uttarkhand,  India
}
\email{swamifma@iitr.ernet.in, mathswami@gmail.com}

\bigskip

\begin{abstract}
Let $\mathcal{W}_{\beta}^\delta(\alpha,\gamma)$ be the class of
normalized analytic functions $f$ defined in the region $|z|<1$
and satisfying

\begin{align*}
{\rm Re\,} e^{i\phi}\left(\dfrac{}{}(1\!-\!\alpha\!+\!2\gamma)\!\left({f}/{z}\right)^\delta
+\left(\alpha\!-\!3\gamma+\gamma\left[\dfrac{}{}\left(1-{1}/{\delta}\right)\left({zf'}/{f}\right)+
{1}/{\delta}\left(1+{zf''}/{f'}\right)\right]\right)\right.\\
\left.\dfrac{}{}\left({f}/{z}\right)^\delta \!\left({zf'}/{f}\right)-\beta\right)>0,
\end{align*}
with the conditions $\alpha\geq 0$, $\beta<1$, $\gamma\geq 0$,
$\delta>0$ and $\phi\in\mathbb{R}$. For a non-negative and
real-valued integrable function $\lambda(t)$ with
$\int_0^1\lambda(t) dt=1$, the generalized non-linear
integral transform is defined as
\begin{align*}
V_{\lambda}^\delta(f)(z)= \left(\int_0^1 \lambda(t) \left({f(tz)}/{t}\right)^\delta dt\right)^{1/\delta}.
\end{align*}
The main aim of the present work is to find conditions on the related
parameters such that
$V_\lambda^\delta(f)(z)\in\mathcal{W}_{\beta_1}^{\delta_1}(\alpha_1,\gamma_1)$,
whenever $f\in\mathcal{W}_{\beta_2}^{\delta_2}(\alpha_2,\gamma_2)$. Further, several
interesting applications for specific choices of $\lambda(t)$
are discussed.
\end{abstract}



\maketitle

\pagestyle{myheadings} \markboth{ Satwanti Devi and A. Swaminathan }{
Inclusion properties of Generalized Integral Transform using Duality Techniques
}

\section{introduction}
Let $\mathcal{A}$ be the class of all normalized and
analytic functions $f$ defined in the open unit disk $\mathbb{D}=\{ z\in\mathbb{C}: |z|<1\}$ that have the Taylor series representation
\begin{align*}
f(z)=z+\sum_{n=2}^\infty a_n z^n.
\end{align*}
Let $\lambda(t):[0,1]\rightarrow \mathbb{R}$ be a non-negative
integrable function which satisfies the condition $\int_0^1\lambda(t)
dt=1$. For $f\in\mathcal{A}$,  consider the generalized integral
transform defined by
\begin{align}\label{eq-weighted-integralOperator}
F_\delta(f)(z):=V_{\lambda}^\delta(f)(z)=\left(\int_0^1 \lambda(t) \left(\dfrac{f(tz)}{t}\right)^\delta dt\right)^{1/\delta},
\quad \delta>0\quad {\rm and}\quad z\in\mathbb{D}.
\end{align}
The power appearing in $\eqref{eq-weighted-integralOperator}$ and elsewhere in this manuscript are meant as principal values.
We are interested in the following problem.
\begin{problem}\label{problem-genl-univ}
Find a class of admissible functions $f \in {\mathcal{A}}$ that are carried by the integral operator defined by
$\eqref{eq-weighted-integralOperator}$ to a class of analytic functions.
\end{problem}
Consideration of such problem for various type of non-linear integral operators exist in the literature. For example, the
integral operator
\begin{align*}
I_{g,c,\alpha}\left[\dfrac{}{}f(z)\right] =
\left( \dfrac{c+\alpha}{g^{c}(z)}\int_0^z f^{\alpha}(t) g^{c-1}(t) g'(t) dt\right)^{1/\alpha}
\end{align*}
where $g \in {\mathcal{A}}$, $g(0)=0$, $g'(0)\neq 0$ and $g(z)\neq 0$ for $z\in {\mathbb{D}}\backslash\{0\}$,
was considered by T. Bulboac\u a in \cite[P.58]{Bulboaca} for analyzing various inclusion properties involving interesting
classes of analytic functions. Choosing $g(z)=z$ gives
\begin{align}\label{eq-nonliear-miller-operator}
I_{c,\alpha}\left[\dfrac{}{}f(z)\right] =
\left( \dfrac{c+\alpha}{z^{c}}\int_0^z f^{\alpha}(t) t^{c-1} dt\right)^{1/\alpha}.
\end{align}
This operator has a rich literature and was considered among several authors by S. S.Miller and P.T. Mocanu \cite[P.319]{Miller-operator} (see also
\cite[P.228]{Miller-Mocanu-book}) for studying various inclusion properties. Even though particular values of the operator given by
$\eqref{eq-nonliear-miller-operator}$ may be related to the operator $\eqref{eq-weighted-integralOperator}$ under some restrictions, these
two operators are entirely different. But the existing literature related to the operator $\eqref{eq-nonliear-miller-operator}$ motivate us
to consider suitable classes of analytic functions that can be studied with reference to the operator given by $\eqref{eq-weighted-integralOperator}$.
For this purpose, we define the class $\mathcal{W}_{\beta}^\delta(\alpha,\gamma)$ in the following way.
{\small{
\begin{align*}
\mathcal{W}_{\beta}^\delta(\alpha,\gamma)\!:\!=\!\left\{\!f\!\in\!\mathcal{A}:
{\rm Re\,} e^{i\phi}\left((1\!-\!\alpha\!+\!2\gamma)\!\left(\frac{f}{z}\right)^\delta
\!+\!\left(\alpha\!-\!3\gamma\!+\!\gamma\left[\left(1\!-\!\frac{1}{\delta}\right)\left(\frac{zf'}{f}\right)\!+\!
\frac{1}{\delta}\left(1\!+\!\frac{zf''}{f'}\right)\right]\right)\right.\right.\\
\left.\left.\left(\frac{f}{z}\right)^\delta \!\left(\frac{zf'}{f}\right)-\beta\right)>0,
\,\,z\in\mathbb{D},\,\,\phi\in\mathbb{R}\right\}.
\end{align*}}}
\hspace*{-2.5mm} Note that, for the particular case
$\mathcal{W}_{\beta}^\delta(\alpha,0)\equiv P_\alpha(\delta,\beta)$,  the operator $\eqref{eq-weighted-integralOperator}$
was examined by R. Aghalary et al. \cite{AghalaryUnivalent} using duality techniques.
Further, for the case $\delta=1$, this operator \eqref{eq-weighted-integralOperator} reduces to the one introduced by R. Fournier and
S. Ruscheweyh \cite{FourRusExtremal}, that contains some of the well-known operator such as Bernardi, Komatu and Hohlov as its
special cases for particular choices of $\lambda(t)$, which has been extensively studied by various authors (for details see
\cite{AbeerS, MahnazC, Ali, saigo, DeviPascu, DeviPascuOrder} and references therein).

It is also interesting to note that results related to other particular cases namely,
$\mathcal{W}_{\beta}^1(\alpha,\gamma)\equiv\mathcal{W}_{\beta}(\alpha,\gamma)$, considered by R.M. Ali et al \cite{AbeerS}
and $\mathcal{W}_{\beta}^\delta(\alpha+\delta+\delta\alpha,\delta\alpha)\equiv
R_\alpha(\delta,\beta)$ exist in the literature.
Note that $F:=z\left({f}/{z}\right)^\delta\in R_\alpha(\delta,\beta)\Longleftrightarrow zF'(z)\in P_\alpha(\delta,\beta)$. Hence the class
$R_\alpha(\delta,\beta)$ is closely related to the class $P_\alpha(\delta,\beta)$.
Inclusion properties of $V_\lambda^\delta(f)(z)\in\mathcal{W}_{\beta}^{\delta}(\alpha,\gamma)$ for
subclasses of analytic functions that have geometrical meaning were considered by the authors of this work in \cite{DeviGenlStar, DeviGenlConvex}.
Since it is challenging to solve Problem $\ref{problem-genl-univ}$ completely, in this manuscript the following particular cases are addressed.
\begin{problem}\label{problem-genl-univ-2}
For given $\xi<1$, to obtain the sharp bounds for $\beta$ such that
\begin{enumerate}[{\rm(i)}]
\item $V_\lambda^\delta(f)(z)\in\mathcal{W}_\xi^\delta(1,0)$  whenever $f(z)\in\mathcal{W}_\beta^\delta(\alpha,\gamma)$ and
\item $V_\lambda^\delta(f)(z)\in\mathcal{W}_\xi^\delta(\alpha,\gamma)$ whenever $f(z)\in \mathcal{W}_\beta^\delta(\alpha,\gamma)$.
\end{enumerate}
\end{problem}
Main results related to Problem $\ref{problem-genl-univ-2}$ and its consequences are given in Section $\ref{sec-genl-univ-mainresults}$ whereas the
required proofs are given in Section $\ref{sec-genl-univ-proofs}$ separately to provide the readers a collective view of the results.

\section{Main results and their consequences}\label{sec-genl-univ-mainresults}
The following Lemma is used in obtaining the main results.
\begin{lemma}\label{Lemma_Univalence}{\rm\cite{Stankiewicz}}
Let $\beta_1,\beta_2<1$ and $\phi\in\mathbb{R}$. Then for
the analytic functions $p$ and $q$ defined in the region
$\mathbb{D}$ with $p(0)=1=q(0)$, along with the conditions
${\rm Re}{\,}e^{i\phi}(p(z)-\beta_1)>0$ and ${\rm
Re}(q(z)-\beta_2)>0$ implies ${\rm
Re}{\,}e^{i\phi}\left({\,}(p(z)\ast q(z)){\,}-\beta\right)>0$,
where $(1-\beta)=2(1-\beta_1)(1-\beta_2)$.
\end{lemma}
Here $\ast$ denotes the convolution or Hadamard product of two
normalized analytic functions $f_1(z)=z+\sum_{n=2}^\infty a_n
z^n$ and $f_2(z)=z+\sum_{n=2}^\infty b_n z^n$, defined in $\mathbb{D}$, given by
\begin{align*}
(f_1\ast f_2)(z)=z+\displaystyle\sum_{n=2}^{\infty}a_nb_nz^n.
\end{align*}

The parameters $\mu,\,\nu\geq0$ introduced in \cite{AbeerS} are
used for further analysis that are defined by the following
relations
\begin{align}\label{eq-mu+nu}
\mu\nu=\gamma\quad  \text{and}\quad \mu+\nu=\alpha-\gamma.
\end{align}
Clearly $\eqref{eq-mu+nu}$ leads to two cases.
\begin{itemize}
\item[{\rm{(i)}}] $\gamma=0 \, \Longrightarrow \mu=0,\,
    \nu=\alpha \geq 0$.
\item[{\rm{(ii)}}] $\gamma>0 \, \Longrightarrow \mu>0,\,
    \nu>0$.
\end{itemize}

Following theorem addresses the first question of Problem $\ref{problem-genl-univ-2}$ and the proof of the same is given in Section
$\ref{sec-genl-univ-proofs}$.
\begin{theorem}\label{theorem_Univalence_1}
Let $\gamma\geq0$ $(\mu\geq0,\nu\geq0)$ and $\delta>0$. Further
let $\xi<1$ and $\beta<1$ is defined by the relation
{\small{
\begin{align}\label{Beta-Cond-Generalized:univalence}
\beta\!=\!\left\{\!\!
               \begin{array}{ll}
                1\!-\!\dfrac{1\!-\!\xi}{2}\left(1-\displaystyle\int_0^1\lambda(t)\left[\dfrac{1}{\nu}\int_0^1\dfrac{ds}{1\!+\!ts^{\mu/\delta}}
                +\left(1\!-\!\dfrac{1}{\nu}\right)\int_0^1\int_0^1\dfrac{drds}{1\!+\!tr^{\nu/\delta}s^{\mu/\delta}}\right]dt\right)^{-1} &\gamma>0, \\\\
                \quad\quad1-\dfrac{1\!-\!\xi}{2}\left(1-\displaystyle\int_0^1\lambda(t)\left[\dfrac{1}{\alpha(1+t)}
                +\left(1\!-\!\dfrac{1}{\alpha}\right)\int_0^1\dfrac{dr}{1\!+\!tr^{\alpha/\delta}}\right]dt\right)^{-1} &\gamma=0.
               \end{array}
             \right.
\end{align}}}
Then for $f(z)\in\mathcal{W}_\beta^\delta(\alpha,\gamma)$, the
function
$V_\lambda^\delta(f)(z)\in\mathcal{W}_\xi^\delta(1,0)$. The
value of $\beta$ is sharp.
\end{theorem}
\begin{remark}
\begin{enumerate}[1.]
  \item For $\gamma=0$, {\rm Theorem
      \ref{theorem_Univalence_1}} reduces to {\rm
      \cite[Theorem 2.1]{AghalaryUnivalent}}.
  \item For $\delta=1$, {\rm Theorem
      \ref{theorem_Univalence_1}} gives the result of {\rm
      \cite[Theorem 2.1]{SarikaUnivalent}}.
\end{enumerate}
\end{remark}
The integral operator \eqref{eq-weighted-integralOperator}
defined by the weight function
\begin{align}\label{eq-bernardi-lambda}
\lambda(t)=(1+c)t^c,\quad c>-1,
\end{align}
is known as generalized Bernardi operator denoted by
$(\mathcal{B}_c^\delta)$. This operator is the particular case
of the generalized integral operators, considered in the work
of R. Aghalary et al. \cite{AghalaryUnivalent} that follows in
the sequel. The operator corresponding to the value $\delta=1$,
was introduced by S. D. Bernardi \cite{bernardi}. Now, using
the integral operator $\mathcal{B}_c^\delta$, the following
corollary is stated as under.
\begin{corollary}
Let $\gamma\geq0$ $(\mu\geq0,\nu\geq0)$, $\delta>0$ and $c>-1$.
Further let $\xi<1$ and $\beta<1$ is defined by the relation

{\footnotesize{
\begin{align*}
\beta=\left\{\!\!
               \begin{array}{ll}
                1-\dfrac{(1-\xi)(2+c)(\delta+\mu)}{2\delta(1+c)}\left[\dfrac{1}{\nu}  {\,}_{3}F_2\left(\!\!\!\!
\begin{array}{cll}&\displaystyle 1,(2+c),\left(1+\dfrac{\delta}{\mu}\right)\\
&\displaystyle (3+c),\left(2+\dfrac{\delta}{\mu}\right)
\end{array};-1\right)\right.\\
\quad\quad\quad\left.+\dfrac{\delta}{(\delta+\nu)} \left(1-\dfrac{1}{\nu}\right){\,}_{4}F_3\left(\!\!\!\!
\begin{array}{cll}&\displaystyle 1,(2+c),\left(1+\dfrac{\delta}{\mu}\right),\left(1+\dfrac{\delta}{\nu}\right)
\\
&\displaystyle (3+c),\left(2+\dfrac{\delta}{\mu}\right),\left(2+\dfrac{\delta}{\nu}\right)
\end{array};-1\right)
\right]^{-1}, &\gamma>0, \\\\
                \quad\quad 1-\dfrac{(1-\xi)(2+c)}{2(1+c)}\left[\dfrac{1}{\alpha}  {\,}_{2}F_1\left(\!\!\!\!
\begin{array}{cll}&\displaystyle 1,(2+c)\\
&\displaystyle (3+c)
\end{array};-1\right)\right.\\
\quad\quad\quad\quad\quad\left.+\dfrac{\delta}{(\delta+\alpha)} \left(1-\dfrac{1}{\alpha}\right){\,}_{3}F_2\left(\!\!\!\!
\begin{array}{cll}&\displaystyle 1,(2+c),\left(1+\dfrac{\delta}{\alpha}\right)\\
&\displaystyle (3+c),\left(2+\dfrac{\delta}{\alpha}\right)
\end{array};-1\right)
\right]^{-1}, &\gamma=0.
               \end{array}
             \right.
\end{align*}}}
Then for $f(z)\in\mathcal{W}_\beta^\delta(\alpha,\gamma)$, the
function
\begin{align*}
\mathcal{B}_c^\delta(f)(z)=\left(\dfrac{(1+c)}{z^{c-\delta+1}}\int_0^z \omega^{c-\delta}(f(\omega))^\delta d\omega\right)^{1/\delta},
\end{align*}
belongs to $\mathcal{W}_\xi^\delta(1,0)$. The value of $\beta$
is sharp.
\end{corollary}
Here $_{p}F_{q}(c_1,\ldots,c_p;d_1,\ldots,d_q;z)$ or $_{p}F_{q}$ denotes the generalized hypergeometric function given by
\begin{align}\label{generalized-hyper-series}
{\,}_{p}F_q\left(\!\!\!\!\!
\begin{array}{cll}&\displaystyle c_1,\ldots,c_p
\\
&\displaystyle d_1,\ldots,d_q
\end{array};z\right)
=\sum_{n=0}^\infty
\dfrac{(c_1)_n\ldots,(c_p)_n}{(d_1)_n\ldots,(d_q)_nn!}z^n,\quad z\in\mathbb{D},
\end{align}
where  $c_i$ $(i=0,1, \ldots,p)$ and $d_j$ $(j=0,1,\ldots,q)$ are the complex parameters with
$d_j\neq0,-1,\ldots$ and $p\leq q+1$. In particular, $_{2}F_{1}$ is the well-known Gaussian
hypergeometric function. For any natural number $n$, the Pochhammer symbol or shifted factorial $(\varepsilon)_n$ is
defined as $(\varepsilon)_0=1$ and $(\varepsilon)_n=\varepsilon(\varepsilon+1)_{n-1}$. The generalized hypergeometric series $_{p}F_q$  defined in
\eqref{generalized-hyper-series}, converges absolutely for all $z$ in $|z|<\infty$ if $p<q+1$, and for $z\in\mathbb{D}$ if $p=q+1$.

The following result gives conditions such that
$V_\lambda^\delta(f)(z)\in\mathcal{W}_\xi^\delta(\alpha,\gamma)$ whenever $f(z)\in \mathcal{W}_\beta^\delta(\alpha,\gamma)$ which addresses
the second part of Problem $\ref{problem-genl-univ-2}$.
\begin{theorem}\label{theorem_Univalence_2}
Let $\gamma\geq0$ $(\mu\geq0,\nu\geq0)$ and $\delta>0$. Further
let $\xi<1$ and $\beta<1$ is defined by the relation
\begin{align}\label{eq-generalized-betacond-main2}
\dfrac{\beta}{(1-\beta)}=-\int_0^1\lambda(t)\dfrac{\left(1-\frac{1+\xi}{1-\xi}t\right)}{1+t}dt.
\end{align}
Then for $f(z)\in\mathcal{W}_\beta^\delta(\alpha,\gamma)$, the
function
$V_\lambda^\delta(f)(z)\in\mathcal{W}_\xi^\delta(\alpha,\gamma)$.
The value of $\beta$ is sharp.
\end{theorem}
It is interesting to note that Theorem $\ref{theorem_Univalence_2}$ cannot be reduced to Theorem $\ref{theorem_Univalence_1}$ for
particular values of the parameters. However, for particular values of the parameters given in Theorem $\ref{theorem_Univalence_2}$ several
results exist in the literature and few of them are listed below.
\begin{remark}
\begin{enumerate}[1.]
\item For $\gamma=0$, $\delta\!=\!1$, $\alpha\!=\!1$ {\rm
    Theorem \ref{theorem_Univalence_2}} reduces to
    {\rm\cite[Theorem 2]{FourRusExtremal}}.
\item For $\gamma=0$ and $\delta=1$, {\rm Theorem
    \ref{theorem_Univalence_2}} gives the result of
    {\rm\cite[Theorem 2.6]{KimRonning}}.
\item For $\gamma=0$, {\rm Theorem
    \ref{theorem_Univalence_2}} reduces to
    {\rm\cite[Theorem 2.2]{AghalaryUnivalent}}.
\item For $\delta=1$, {\rm Theorem
    \ref{theorem_Univalence_2}} gives the result of
    {\rm\cite[Theorem 2.3]{SarikaUnivalent}}.
\end{enumerate}
\end{remark}
On setting $\lambda(t)$ given in \eqref{eq-bernardi-lambda},
Theorem \ref{theorem_Univalence_2} will lead to the
following result.
\begin{corollary}
Let $\gamma\geq0$ $(\mu\geq0,\nu\geq0)$, $\delta>0$ and $c>-1$.
Further let $\xi<1$ and $\beta<1$ is defined by the relation
\begin{align*}
\beta=\dfrac{2(1+c)\,_2F_1(1,2+c;3+c;-1)-(2+c)(1-\xi)}{2(1+c)\,_2F_1(1,2+c;3+c;-1)}.
\end{align*}
Then for $f(z)\in\mathcal{W}_\beta^\delta(\alpha,\gamma)$, the
function
$\mathcal{B}_c^\delta(f)(z)\in\mathcal{W}_\xi^\delta(\alpha,\gamma)$.
The value of $\beta$ is sharp.
\end{corollary}
Let
\begin{align*}
\lambda(t)=\dfrac{\Gamma(c)}{\Gamma(a)\Gamma(b)\Gamma(c-a-b+1)}
t^{b-1}(1-t)^{c-a-b}{\,}_{2}F_1\left(\!\!\!\!
\begin{array}{cll}&\displaystyle c-a,\quad 1-a
\\
&\displaystyle c-a-b+1
\end{array};1-t\right),
\end{align*}
then the integral operator \eqref{eq-weighted-integralOperator}
defined by the above weight function $\lambda(t)$ is the known
as generalized Hohlov operator denoted by
$\mathcal{H}_{a,b,c}^\delta$. This integral operator was
considered in the work of A. Ebadian et al. \cite{Aghalary}. Its
representation in the form of convolution is given as
$\mathcal{H}:=\mathcal{H}_{a,b,c}^\delta(f)(z)=(z^\delta{\,}_2F_1(a,b;c;z))\ast
(f(z))^\delta$. For the case $\delta=1$, the convolution form
of the reduced integral transform is given in the work of Y. C.
Kim and F. Ronning \cite{KimRonning} and studied by several
authors later. The operator $\mathcal{H}_{a,b,c}^\delta$ with
$a=1$ is the generalized Carlson-Shaffer operator
($\mathcal{L}_{b,c}^\delta$) \cite{CarlsonShaffer}.
Corresponding to these operators the following results are
obtained.
\begin{theorem}\label{thm-univ:hyper}
Let $\alpha>\gamma\left(1+\frac{(2a+1)}{\delta}\right)\geq0$,
$\delta>0$,
$0<a\leq\min\left\{1,\frac{\delta}{2(\alpha-\gamma)}\right\}$
and $0<1+b<(c-a)<2$. Then for
$f(z)\in\mathcal{W}_{\beta_1}^\delta(0,0)$, the function
$\mathcal{H}_{a,b,c}^\delta(f)(z)$ belongs to the class
$\mathcal{W}_\beta^\delta(\alpha,\gamma)$, where
$\beta=1-2(1-\beta_1)(1-\beta_2)$ with
{\small{
\begin{align*}
\beta_2:=\left[1\!-\!\dfrac{a}{\delta}\left(\alpha\!-\!\gamma\left(1\!+\!\dfrac{a}{\delta}\right)\right)\right]\,_2F_1(a,b;c;-1)+
\dfrac{a}{\delta}\left[\alpha\!-\!\gamma\left(1\!+\!\dfrac{(2a+1)}{\delta}\right)\right]\,_2F_1(a+1,b;c;-1)\\
+\dfrac{a(a+1)\gamma}{\delta^2}\,_2F_1(a+2,b;c;-1).
\end{align*}}}
The value of $\beta$ is sharp.
\end{theorem}
For $a=1$, Theorem \ref{thm-univ:hyper} lead to the following
result.
\begin{corollary}
Let $\alpha>\gamma\left(1+\frac{3}{\delta}\right)\geq0$,
$\delta>0$, and $0<b<(c-1)<1$. Then for
$f(z)\in\mathcal{W}_{\beta_1}^\delta(0,0)$, the function
$\mathcal{L}_{b,c}^\delta(f)(z)$ belongs to the class
$\mathcal{W}_\beta^\delta(\alpha,\gamma)$, where
$\beta=1-2(1-\beta_1)(1-\beta_2)$ with

\begin{align*}
\beta_2:=\left[1-\dfrac{1}{\delta}\left(\alpha-\gamma\left(1+\dfrac{1}{\delta}\right)\right)\right]\,_2F_1(1,b;c;-1)+
\dfrac{1}{\delta}\left[\alpha-\gamma\left(1+\dfrac{3}{\delta}\right)\right]\,_2F_1(2,b;c;-1)\\
+\dfrac{2\gamma}{\delta^2}\,_2F_1(3,b;c;-1).
\end{align*}
The value of $\beta$ is sharp.
\end{corollary}

It is important to note that similar results for various operators considered in \cite{DeviGenlStar, DeviGenlConvex} by the authors of
this work are already obtained and are appearing elsewhere.

\section{Proofs of theorems $\ref{theorem_Univalence_1}$, $\ref{theorem_Univalence_2}$ and $\ref{thm-univ:hyper}$}\label{sec-genl-univ-proofs}
\noindent
{\underline{\bf{Proof of Theorem $\ref{theorem_Univalence_1}$}}}.
Since the case $\gamma=0(\mu=0, \nu=\alpha>0)$ corresponds to
\cite[Theorem 2.1]{AghalaryUnivalent}, it is enough to obtain
the condition for $\gamma>0$. Let
{\small{
\begin{align}\label{eq-generalized:starlike-H}
H(z):=(1\!-\!\alpha\!+\!2\gamma)\!\left(\frac{f}{z}\right)^\delta
\!+\!\left(\alpha\!-\!3\gamma+\gamma\left[\left(1\!-\!\frac{1}{\delta}\right)\left(\frac{zf'}{f}\right)\!+\!
\frac{1}{\delta}\left(1\!+\!\frac{zf''}{f'}\right)\right]\right)
\left(\frac{f}{z}\right)^\delta\!\left(\frac{zf'}{f}\right).
\end{align}}}
Now, considering
$ \displaystyle
P(z)=\left(\dfrac{f}{z}\right)^\delta
$
a simple computation gives
%
%
\begin{align}\label{eq-P-form1}
\left(1-\dfrac{1}{\delta}\right)P(z)+\left(\dfrac{1}{\delta}\right)(zP(z))'=\left(\dfrac{f}{z}\right)^\delta\left(\dfrac{zf'}{f}\right).
\end{align}
%
Now, differentiating \eqref{eq-P-form1} twice and applying successively in
%
%
\eqref{eq-generalized:starlike-H} leads
to
\begin{align*}
H(z)=\left(1-\dfrac{(\alpha-\gamma)}{\delta}+\dfrac{\gamma}{\delta^2}\right)P(z)+
\left(\dfrac{(\alpha-\gamma)}{\delta}-\dfrac{2\gamma}{\delta^2}\right)(zP(z))'+\dfrac{\gamma}{\delta^2}(z(zP(z))')'
\end{align*}
which upon using \eqref{eq-mu+nu} leads to
\begin{align*}
H(z)=\left(1-\dfrac{(\mu+\nu)}{\delta}+\dfrac{\mu\nu}{\delta^2}\right)P(z)+
\left(\dfrac{(\mu+\nu)}{\delta}-\dfrac{2\mu\nu}{\delta^2}\right)(zP(z))'+\dfrac{\mu\nu}{\delta^2}(z(zP(z))')'.
\end{align*}
Considering a series expansion of the form $P(z)=1+\sum_{n=1}^\infty b_nz^n$ in the above equality gives
\begin{align}\label{eq-generalized-H(z)-seriesform}
H(z)=1+\dfrac{1}{\delta^2}\sum_{n=1}^\infty (\delta+n\mu)(\delta+n\nu)b_nz^n.
\end{align}
From \eqref{eq-P-form1} and \eqref{eq-generalized-H(z)-seriesform}, it is a simple exercise
to see that
\begin{align}\label{eq-zf'/f-convol-form}
\left(\dfrac{f}{z}\right)^\delta\left(\dfrac{zf'}{f}\right)=\left(1+\delta\sum_{n=1}^\infty
\dfrac{(n+\delta)}{(\delta+n\mu)(\delta+n\nu)} z^n\right)
\ast H(z).
\end{align}
Since
\begin{align*}
1+\delta\sum_{n=1}^\infty\dfrac{(n+\delta)}{(\delta+n\mu)(\delta+n\nu)} z^n=&
\left[\left(\dfrac{1}{\delta}\right){\,}_3F_2\left(2,\dfrac{\delta}{\nu},\dfrac{\delta}{\mu};
\left(1+\dfrac{\delta}{\nu}\right),\left(1+\dfrac{\delta}{\mu}\right);z\right)\right.\\
&\left.+\left(1-\dfrac{1}{\delta}\right){\,}_3F_2\left(1,\dfrac{\delta}{\nu},\dfrac{\delta}{\mu};
\left(1+\dfrac{\delta}{\nu}\right),\left(1+\dfrac{\delta}{\mu}\right);z\right)\right],
\end{align*}
\eqref{eq-zf'/f-convol-form} is equivalent to
\begin{align}\label{eq-zf'/f-hypergeometric-form}
\left(\dfrac{f}{z}\right)^\delta\left(\dfrac{zf'}{f}\right)=&\left[\left(\dfrac{1}{\delta}\right){\,}_3F_2\left(2,\dfrac{\delta}{\nu},\dfrac{\delta}{\mu};
\left(1+\dfrac{\delta}{\nu}\right),\left(1+\dfrac{\delta}{\mu}\right);z\right)\right.\nonumber\\
&\left.+\left(1-\dfrac{1}{\delta}\right){\,}_3F_2\left(1,\dfrac{\delta}{\nu},\dfrac{\delta}{\mu};
\left(1+\dfrac{\delta}{\nu}\right),\left(1+\dfrac{\delta}{\mu}\right);z\right)\right]
\ast H(z).
\end{align}
Taking the logarithmic derivative on both sides of \eqref{eq-weighted-integralOperator} we have

\begin{align}\label{eq-generalized-F-integral-rep-1}
\dfrac{zF_\delta'(z)}{F_\delta(z)}\left(\dfrac{F_\delta(z)}{z}\right)^\delta
=\int_0^1\dfrac{\lambda(t)}{1-tz}dt \ast \left(\left(\dfrac{zf'(z)}{f(z)}\right)\left(\dfrac{f(z)}{z}\right)^{\delta}\right).
\end{align}
Now, substituting \eqref{eq-zf'/f-hypergeometric-form} in
\eqref{eq-generalized-F-integral-rep-1} will give

\begin{align*}
\dfrac{zF_\delta'(z)}{F_\delta(z)}\left(\dfrac{F_\delta(z)}{z}\right)^\delta
=\int_0^1\dfrac{\lambda(t)}{1-tz}dt \ast \left[\left(\dfrac{1}{\delta}\right){\,}_3F_2\left(2,\dfrac{\delta}{\nu},\dfrac{\delta}{\mu};
\left(1+\dfrac{\delta}{\nu}\right),\left(1+\dfrac{\delta}{\mu}\right);z\right)\right.\nonumber\\
\left.+\left(1-\dfrac{1}{\delta}\right){\,}_3F_2\left(1,\dfrac{\delta}{\nu},\dfrac{\delta}{\mu};
\left(1+\dfrac{\delta}{\nu}\right),\left(1+\dfrac{\delta}{\mu}\right);z\right)\right]\ast H(z).
\end{align*}
For $f(z)\in\mathcal{W}_\beta^\delta(\alpha,\gamma)$, then it
is easy to see that for some $\phi\in\mathbb{R}$, ${\rm
Re}{\,}e^{i\phi}(H(z)-\beta)>0$. Therefore, for $\gamma>0$, it
is required to prove the claim that
\begin{align}\label{eq-generalized-hyper-main1}
{\rm Re}\int_0^1\lambda(t)&\left[\left(\dfrac{1}{\delta}\right){\,}_3F_2\left(2,\dfrac{\delta}{\nu},\dfrac{\delta}{\mu};
\left(1+\dfrac{\delta}{\nu}\right),\left(1+\dfrac{\delta}{\mu}\right);tz\right)\right.\nonumber\\
&\left.+\left(1-\dfrac{1}{\delta}\right){\,}_3F_2\left(1,\dfrac{\delta}{\nu},\dfrac{\delta}{\mu};
\left(1+\dfrac{\delta}{\nu}\right),\left(1+\dfrac{\delta}{\mu}\right);tz\right)\right]dt>1-\dfrac{1-\xi}{2(1-\beta)},
\end{align}
which by applying Lemma \ref{Lemma_Univalence}, implies that
$F_\delta\in\mathcal{W}_\xi^\delta(1,0)$. Now it is enough to
verify inequality \eqref{eq-generalized-hyper-main1}.

From the identity
\begin{align*}
_3F_2(2,a,b;c,d;z)=(c-1){\,}_3F_2(1,a,b;c-1,d;z)-(c-2){\,}_3F_2(1,a,b;c,d;z),
\end{align*}
it follows that
\begin{align*}
_3F_2\left(2,\dfrac{\delta}{\nu},\dfrac{\delta}{\mu};
\left(1+\dfrac{\delta}{\nu}\right),\left(1+\dfrac{\delta}{\mu}\right);z\right)
=\left(\dfrac{\delta}{\nu}\right){\,}_2F_1\left(1,\dfrac{\delta}{\mu};\left(1+\dfrac{\delta}{\mu}\right);z\right)\\
-\left(\dfrac{\delta}{\nu}-1\right){\,}_3F_2\left(1,\dfrac{\delta}{\nu},\dfrac{\delta}{\mu};
\left(1+\dfrac{\delta}{\nu}\right),\left(1+\dfrac{\delta}{\mu}\right);z\right).
\end{align*}
Therefore
{\footnotesize{
\begin{align*}
{\rm Re}\int_0^1\lambda(t)&\left[\left(\dfrac{1}{\delta}\right){\,}_3F_2\left(2,\dfrac{\delta}{\nu},\dfrac{\delta}{\mu};
\left(1+\dfrac{\delta}{\nu}\right),\left(1+\dfrac{\delta}{\mu}\right);tz\right)\right.\\
&\quad\quad\left.+\left(1-\dfrac{1}{\delta}\right){\,}_3F_2\left(1,\dfrac{\delta}{\nu},\dfrac{\delta}{\mu};
\left(1+\dfrac{\delta}{\nu}\right),\left(1+\dfrac{\delta}{\mu}\right);tz\right)\right]dt\\
={\rm Re}\int_0^1\!\!\lambda(t)&\left[\left(\dfrac{1}{\nu}\right)\,_2F_1\left(1,\dfrac{\delta}{\mu};
\left(1+\dfrac{\delta}{\mu}\right);tz\right)+\left(1-\dfrac{1}{\nu}\right)\,_3F_2\left(1,\dfrac{\delta}{\nu},\dfrac{\delta}{\mu};
\left(1+\dfrac{\delta}{\nu}\right),\left(1+\dfrac{\delta}{\mu}\right);tz\right)\right]dt.
\end{align*}}}
Since the integral form of the following generalized
hypergeometric function is given as
\begin{align*}
_2F_1\left(1,\dfrac{\delta}{\mu};\left(1+\dfrac{\delta}{\mu}\right);tz\right)=\int_0^1\dfrac{ds}{(1-tzs^{\mu/\delta})}
\end{align*}
and
\begin{align*}
_3F_2\left(1,\dfrac{\delta}{\nu},\dfrac{\delta}{\mu};\left(1+\dfrac{\delta}{\nu}\right),
\left(1+\dfrac{\delta}{\mu}\right);tz\right)=\int_0^1\int_0^1\dfrac{drds}{(1-tzr^{\nu/\delta}s^{\mu/\delta})},
\end{align*}
therefore
{\footnotesize{
\begin{align*}
{\rm Re}\int_0^1\lambda(t)\left[\left(\dfrac{1}{\nu}\right)\,_2F_1\left(1,\dfrac{\delta}{\mu};
\left(1+\dfrac{\delta}{\mu}\right);tz\right)+\left(1-\dfrac{1}{\nu}\right)\,_3F_2\left(1,\dfrac{\delta}{\nu},\dfrac{\delta}{\mu};
\left(1+\dfrac{\delta}{\nu}\right),\left(1+\dfrac{\delta}{\mu}\right);tz\right)\right]dt\\
={\rm Re}\int_0^1\lambda(t)\left[\left(\dfrac{1}{\nu}\right)\int_0^1\dfrac{ds}{(1-tzs^{\mu/\delta})}
+\left(1-\dfrac{1}{\nu}\right)\int_0^1\int_0^1\dfrac{drds}{(1-tzr^{\nu/\delta}s^{\mu/\delta})}\right]dt.
\end{align*}}}
It is evident that ${\rm Re}
\left(\frac{1}{1-tz}\right)>\frac{1}{(1+t)}$, for $|z|<1$. Thus
by the given hypothesis
\begin{align*}
&{\rm Re}\int_0^1\lambda(t)\left[\left(\dfrac{1}{\nu}\right)\int_0^1\dfrac{ds}{(1-tzs^{\mu/\delta})}
+\left(1-\dfrac{1}{\nu}\right)\int_0^1\int_0^1\dfrac{drds}{(1-tzr^{\nu/\delta}s^{\mu/\delta})}\right]dt\\
\geq&\int_0^1\lambda(t)\left[\left(\dfrac{1}{\nu}\right)\int_0^1\dfrac{ds}{(1+ts^{\mu/\delta})}
+\left(1-\dfrac{1}{\nu}\right)\int_0^1\int_0^1\dfrac{drds}{(1+tr^{\nu/\delta}s^{\mu/\delta})}\right]dt\\
=&1-\dfrac{1-\xi}{2(1-\beta)}
\end{align*}
and the proof is complete.

Now, to verify the sharpness let
$f(z)\in\mathcal{W}_\beta^\delta(\alpha,\gamma)$, therefore it
satisfies the differential equation
\begin{align}\label{eq-generalized-f/z-extremal}
\dfrac{\mu\nu}{\delta^2}{\,} z^{1-\delta/\mu}\left(z^{\delta/\mu-\delta/\nu+1}
\left( z^{\delta/\nu}\left(\dfrac{f}{z}\right)^\delta\right)'\right)'
=\beta+(1-\beta)\dfrac{1+z}{1-z}
\end{align}
with $\beta<1$ defined in
\eqref{Beta-Cond-Generalized:univalence}. From \eqref{eq-mu+nu},
\eqref{eq-generalized-f/z-extremal} and the series representation of $P(z)$, an easy
calculation gives
\begin{align*}
\left(\frac{f(z)}{z}\right)^\delta\left(\frac{zf'(z)}{f(z)}\right)
=1+2(1-\beta)\sum_{n=1}^\infty \dfrac{\delta (n+\delta)z^n}{(\delta+n\nu)(\delta+n\mu)}.
\end{align*}
Therefore from \eqref{eq-generalized-F-integral-rep-1}, we have
\begin{align}\label{eq-zF'/F-F/z-series-form}
\left(\frac{F_\delta(z)}{z}\right)^\delta\left(\frac{zF_\delta'(z)}{F_\delta(z)}\right)
&=\int_0^1\lambda(t)\left(1+2(1-\beta)\sum_{n=1}^\infty \dfrac{\delta (n+\delta)(tz)^n}{(\delta+n\nu)(\delta+n\mu)}\right)dt\nonumber\\
&=1+2(1-\beta)\sum_{n=1}^\infty \dfrac{\delta (n+\delta)\tau_nz^n}{(\delta+n\nu)(\delta+n\mu)}
\end{align}
where $\tau_n=\int_0^1t^n\lambda(t) dt$. By the simple
adjustment, \eqref{Beta-Cond-Generalized:univalence} can be
written as
\begin{align*}
\dfrac{1}{(1-\beta)}=\dfrac{2}{(1-\xi)}
\left(1-\displaystyle\int_0^1\lambda(t)\left[\dfrac{1}{\nu}\int_0^1\dfrac{ds}{1+ts^{\mu/\delta}}
+\left(1-\dfrac{1}{\nu}\right)\int_0^1\int_0^1\dfrac{drds}{1+tr^{\nu/\delta}s^{\mu/\delta}}\right]dt\right)
\end{align*}
or equivalently,
{\small{
\begin{align*}
\dfrac{1}{(1-\beta)}=\dfrac{2}{(1-\xi)}
\left(1-\displaystyle\int_0^1\lambda(t)\left[\dfrac{\delta}{\nu}\sum_{n=0}^\infty\dfrac{(-1)^nt^n}{(\delta+n\mu)}
+\delta^2\left(1-\dfrac{1}{\nu}\right)\sum_{n=0}^\infty\dfrac{(-1)^nt^n}{(\delta+n\mu)(\delta+n\nu)}\right]dt\right).
\end{align*}}}
As $\tau_n=\int_0^1t^n\lambda(t) dt$,
\begin{align*}
\dfrac{1}{(1-\beta)}=\dfrac{2}{(1-\xi)}
\left(1-\dfrac{\delta}{\nu}\sum_{n=0}^\infty\dfrac{(-1)^n\tau_n}{(\delta+n\mu)}
-\delta^2\left(1-\dfrac{1}{\nu}\right)\sum_{n=0}^\infty\dfrac{(-1)^n\tau_n}{(\delta+n\mu)(\delta+n\nu)}\right)
\end{align*}
or
\begin{align}\label{eq-xi-beta:series:relat}
\xi=1+2(1-\beta)\sum_{n=1}^\infty\dfrac{(-1)^n\delta(n+\delta)\tau_n}{(\delta+n\mu)(\delta+n\nu)}.
\end{align}
Further using \eqref{eq-zF'/F-F/z-series-form} and
\eqref{eq-xi-beta:series:relat}, we have
\begin{align*}
\left.\left(\frac{F_\delta(z)}{z}\right)^\delta\left(\frac{zF_\delta'(z)}{F_\delta(z)}\right)\right|_{z=-1}
=1+2(1-\beta)\sum_{n=1}^\infty\dfrac{(-1)^n\delta(n+\delta)\tau_n}{(\delta+n\mu)(\delta+n\nu)}=\xi,
\end{align*}
which clearly implies the sharpness of the result.
\qed

\noindent
{\underline{\bf{Proof of Theorem $\ref{theorem_Univalence_2}$}}}.
Since the case $\gamma=0(\mu=0, \nu=\alpha>0)$ corresponds to
\cite[Theorem 2.2]{AghalaryUnivalent}, therefore the case
$\gamma>0$ is taken into consideration. Now, for the function
$f(z)\in\mathcal{W}_\beta^\delta(\alpha,\gamma)$, let
\begin{align}\label{eq-generalized-G(z)}
G(z)=\dfrac{H(z)-\beta}{(1-\beta)}.
\end{align}
Then it is easy to see that for some $\phi\in\mathbb{R}$, ${\rm
Re}\left(e^{i\phi}{\,}G(z)\right)>0$, where $H(z)$ is defined
in \eqref{eq-generalized:starlike-H}.

Now, the following two cases are discussed.

At first, let $\gamma\neq\delta\alpha/(3\delta-1)$. From
\eqref{eq-generalized:starlike-H},  \eqref{eq-generalized-F-integral-rep-1} and
\eqref{eq-generalized-G(z)},
%
{\small{
\begin{align}\label{eq-generalized-G(z)2}
\left(\dfrac{F_\delta(z)}{z}\right)^\delta \left(\dfrac{zF_\delta'(z)}{F_\delta(z)}\right)=\int_0^1\!\!\!\dfrac{\lambda(t)}{(1-tz)}dt\ast
\dfrac{1}{\left(\alpha\!-\!3\gamma\!+\!\gamma/\delta\right)}
\left(\beta+(1-\beta)G(z)-(1\!-\!\alpha\!+\!2\gamma)\left(\dfrac{f}{z}\right)^\delta\right.\nonumber\\
\left.-\gamma\left(\left(1-\dfrac{1}{\delta}\right)\left(\dfrac{zf'}{f}\right)+
\dfrac{1}{\delta}\left(\dfrac{zf''}{f'}\right)\right)\left(\dfrac{f}{z}\right)^\delta
\left(\dfrac{zf'}{f}\right)\right).
\end{align}}}
Now, from \eqref{eq-weighted-integralOperator} and
\eqref{eq-generalized-F-integral-rep-1}, it is a simple
exercise to see that
{\small{
\begin{align}\label{eq-generalized-F-integral-rep-2}
&\left(\left(1-\dfrac{1}{\delta}\right)\left(\dfrac{zF_\delta'(z)}{F_\delta(z)}\right)+
\dfrac{1}{\delta}\left(\dfrac{zF_\delta''(z)}{F_\delta'(z)}\right)\right)\left(\dfrac{F_\delta(z)}{z}\right)^\delta
\left(\dfrac{zF_\delta'(z)}{F_\delta(z)}\right)\nonumber\\
&\quad\quad=\int_0^1\dfrac{\lambda(t)}{1-tz}dt \ast \left(\left(1-\dfrac{1}{\delta}\right)\left(\dfrac{zf'(z)}{f(z)}\right)+
\dfrac{1}{\delta}\left(\dfrac{zf''(z)}{f'(z)}\right)\right)\left(\dfrac{f(z)}{z}\right)^\delta
\left(\dfrac{zf'(z)}{f(z)}\right).
\end{align}}}
Therefore \eqref{eq-generalized-G(z)2} and
\eqref{eq-generalized-F-integral-rep-2} leads to
{\footnotesize{
\begin{align*}
\left(\dfrac{F_\delta(z)}{z}\right)^\delta \left(\dfrac{zF_\delta'(z)}{F_\delta(z)}\right)=
&\left[\left(\beta+(1-\beta)\int_0^1\dfrac{\lambda(t)}{(1-tz)}dt\right)\ast
\dfrac{1}{\left(\alpha-3\gamma+\gamma/\delta\right)}G(z)\right]\nonumber\\
&-\dfrac{1}{\left(\alpha-3\gamma+\gamma/\delta\right)}
\left((1-\alpha+2\gamma)\left(\dfrac{F_\delta(z)}{z}\right)^\delta
+\gamma\left(\left(1-\dfrac{1}{\delta}\right)\left(\dfrac{zF_\delta'(z)}{F_\delta(z)}\right)\right.\right.\\
&\quad\quad\left.\left.+\dfrac{1}{\delta}\left(\dfrac{zF_\delta''(z)}{F_\delta'(z)}\right)\right)\left(\dfrac{F_\delta(z)}{z}\right)^\delta
\left(\dfrac{zF_\delta'(z)}{F_\delta(z)}\right)\right)
\end{align*}}}
or equivalently,
\newline
$\displaystyle
\left(\beta+(1-\beta)\int_0^1\dfrac{\lambda(t)}{(1-tz)}dt\right)\ast G(z)
=(1\!-\!\alpha\!+\!2\gamma)\left(\frac{F_\delta(z)}{z}\right)^\delta
$
\begin{align}\label{eq-generalized-G(z)3}
+\left(\alpha\!-\!3\gamma\!+\!\gamma\left[\left(1\!-\!\frac{1}{\delta}\right)\left(\frac{zF_\delta'(z)}{F_\delta(z)}\right)+
\frac{1}{\delta}\left(1\!+\!\frac{zF_\delta''(z)}{F_\delta'(z)}\right)\right]\right)
\left(\frac{F_\delta(z)}{z}\right)^\delta\left(\frac{zF_\delta'(z)}{F_\delta(z)}\right).\nonumber\\
\end{align}
Now consider $\gamma=\delta\alpha/(3\delta-1)$. It is easy to
see that the equation \eqref{eq-generalized-G(z)} is equivalent
to
\begin{align*}
\left(1-\gamma+\dfrac{\gamma}{\delta}\right)\left(\dfrac{f}{z}\right)^\delta
+\gamma\left(\left(1-\dfrac{1}{\delta}\right)\left(\dfrac{zf'}{f}\right)+
\dfrac{1}{\delta}\left(\dfrac{zf''}{f'}\right)\right)\left(\dfrac{f}{z}\right)^\delta
\left(\dfrac{zf'}{f}\right)\nonumber\\
=\beta+(1-\beta)G(z).
\end{align*}
Therefore the above expression along with
\eqref{eq-generalized-F-integral-rep-1} and
\eqref{eq-generalized-F-integral-rep-2} gives
\newline
$\displaystyle
G(z)\ast\left(\beta+(1-\beta)\int_0^1\dfrac{\lambda(t)}{(1-tz)}dt\right)
=\left(1-\gamma+\dfrac{\gamma}{\delta}\right)\left(\dfrac{F_\delta(z)}{z}\right)^\delta
$
\begin{align*}
+\gamma\left(\left(1-\dfrac{1}{\delta}\right)\left(\dfrac{zF_\delta'(z)}{F_\delta(z)}\right)+
\dfrac{1}{\delta}\left(\dfrac{zF_\delta''(z)}{F_\delta'(z)}\right)\right)\left(\dfrac{F_\delta(z)}{z}\right)^\delta
\left(\dfrac{zF_\delta'(z)}{F_\delta(z)}\right),
\end{align*}
which coincides with \eqref{eq-generalized-G(z)3}, when
$\gamma=\delta\alpha/(3\delta-1)$.

Moreover, $F_\delta(z)\in\mathcal{W}_\xi^\delta(\alpha,\gamma)$
if, and only if, the functions $F_\delta(z)$ and $J(z)$ are
defined by the relation
\begin{align}\label{eq-generalized-J(z)-F(z)}
\left(\dfrac{J(z)}{z}\right)^\delta=\dfrac{\left(\frac{F_\delta(z)}{z}\right)^\delta-\xi}{(1-\xi)}\quad
\Longleftrightarrow\quad\left(\frac{F_\delta(z)}{z}\right)^\delta=(1-\xi)\left(\dfrac{J(z)}{z}\right)^\delta+\xi.
\end{align}
From the above expression, a simple computation give
\begin{align}\label{eq-generalized-J(z)-F(z)1}
\left(\frac{F_\delta(z)}{z}\right)^\delta\left(\frac{zF_\delta'(z)}{F_\delta(z)}\right)
=(1-\xi)\left(\dfrac{J(z)}{z}\right)^\delta\left(\frac{zJ'(z)}{J(z)}\right)+\xi .
\end{align}
Further using \eqref{eq-generalized-J(z)-F(z)1} will give
\begin{align}\label{eq-generalized-J(z)-F(z)2}
&\left(\left(1-\dfrac{1}{\delta}\right)\left(\dfrac{zF_\delta'(z)}{F_\delta(z)}\right)+
\dfrac{1}{\delta}\left(\dfrac{zF_\delta''(z)}{F_\delta'(z)}\right)\right)\left(\dfrac{F_\delta(z)}{z}\right)^\delta
\left(\dfrac{zF_\delta'(z)}{F_\delta(z)}\right)\nonumber\\
&=(1-\xi)\left(\left(1-\dfrac{1}{\delta}\right)\left(\dfrac{zJ'(z)}{J(z)}\right)+
\dfrac{1}{\delta}\left(\dfrac{zJ''(z)}{J'(z)}\right)\right)\left(\dfrac{J(z)}{z}\right)^\delta
\left(\dfrac{zJ'(z)}{J(z)}\right)+\xi\left(1-\delta\right).
\end{align}
Putting the values from \eqref{eq-generalized-J(z)-F(z)},
\eqref{eq-generalized-J(z)-F(z)1} and
\eqref{eq-generalized-J(z)-F(z)2} in
\eqref{eq-generalized-G(z)}, it follows that
\newline
$\displaystyle
G(z)=(1-\alpha+2\gamma)\left(\frac{J(z)}{z}\right)^\delta
$
{\small{
\begin{align*}
+\left(\alpha-3\gamma+\gamma\left[\left(1-\frac{1}{\delta}\right)\left(\frac{zJ'(z)}{J(z)}\right)+
\frac{1}{\delta}\left(1+\frac{zJ''(z)}{J'(z)}\right)\right]\right)
\left(\frac{J(z)}{z}\right)^\delta\left(\frac{zJ'(z)}{J(z)}\right),
\end{align*}}}
which implies that the function
$J(z)\in\mathcal{W}_0^\delta(\alpha,\gamma)$. Further using
\eqref{eq-generalized-G(z)3} will give
\newline
$
\displaystyle
\left(\dfrac{\beta-\xi}{1-\xi}+\dfrac{1-\beta}{1-\xi}\int_0^1\dfrac{\lambda(t)}{(1-tz)}dt\right)\ast G(z)
=
(1\!-\!\alpha\!+\!2\gamma)\left(\frac{J(z)}{z}\right)^\delta
$
{\footnotesize{
\begin{align*}
+\left(\alpha-3\gamma+\gamma\left[\left(1-\frac{1}{\delta}\right)\left(\frac{zJ'(z)}{J(z)}\right)+
\frac{1}{\delta}\left(1+\frac{zJ''(z)}{J'(z)}\right)\right]\right)
\left(\frac{J(z)}{z}\right)^\delta\left(\frac{zJ'(z)}{J(z)}\right).
\end{align*}
}}
As ${\rm Re}\left(e^{i\phi}{\,}G(z)\right)>0$, from
the result given in \cite[P. 23]{Rus} and by the above
expression, it is easy to note that
{\footnotesize{
\begin{align*}
(1\!-\!\alpha\!+\!2\gamma)\left(\frac{J(z)}{z}\right)^\delta
+\left(\alpha\!-\!3\gamma\!+\!\gamma\left[\left(1\!-\!\frac{1}{\delta}\right)\left(\frac{zJ'(z)}{J(z)}\right)+
\frac{1}{\delta}\left(1\!+\!\frac{zJ''(z)}{J'(z)}\right)\right]\right)
\left(\frac{J(z)}{z}\right)^\delta\left(\frac{zJ'(z)}{J(z)}\right)\neq0,
\end{align*}}}
if, and only if,
\begin{align*}
{\rm Re}\left(\dfrac{\beta-\xi}{1-\xi}+\dfrac{1-\beta}{1-\xi}\int_0^1\dfrac{\lambda(t)}{(1-tz)}dt\right)>\dfrac{1}{2}.
\end{align*}
Now, using ${\rm
Re}\left(\frac{1}{1-tz}\right)>\left(\frac{1}{1+t}\right)$,
$|z|<1$, in the above inequality, it follows that
\begin{align*}
{\rm Re}\left(\dfrac{\beta-\xi}{1-\xi}+\dfrac{1-\beta}{1-\xi}\int_0^1\dfrac{\lambda(t)}{(1-tz)}dt\right)>
\dfrac{\beta-\xi}{1-\xi}+\dfrac{1-\beta}{1-\xi}\int_0^1\dfrac{\lambda(t)}{(1+t)}dt.
\end{align*}
Since the expression \eqref{eq-generalized-betacond-main2} can
be rewritten as
\begin{align*}
\dfrac{\beta-\xi}{1-\beta}+\int_0^1\dfrac{\lambda(t)}{(1+t)}dt=\dfrac{(1-\xi)}{2(1-\beta)},
\end{align*}
the above inequality gives
\begin{align*}
{\rm Re}\left(\dfrac{\beta-\xi}{1-\xi}+\dfrac{1-\beta}{1-\xi}\int_0^1\dfrac{\lambda(t)}{(1-tz)}dt\right)>
\dfrac{\beta-\xi}{1-\xi}+\dfrac{1-\beta}{1-\xi}\int_0^1\dfrac{\lambda(t)}{(1+t)}dt=\dfrac{1}{2}
\end{align*}
and this completes the proof.

Now, to verify the sharpness let $f(z)\in\mathcal{W}_\beta^\delta(\alpha,\gamma)$. Using the series representation of $f(z)$
in \eqref{eq-weighted-integralOperator} will give
\begin{align}\label{eq-F/z_Series}
\left(\frac{F_\delta(z)}{z}\right)^\delta
=1+2(1-\beta)\sum_{n=1}^\infty \dfrac{\delta^2\tau_nz^n}{(\delta+n\nu)(\delta+n\mu)},
\end{align}
where $\tau_n=\int_0^1t^n\lambda(t) dt$.
Now, consider
\newline
$\displaystyle
H_0(z):=(1-\alpha+2\gamma)\left(\frac{F_\delta(z)}{z}\right)^\delta
$
{\small{
\begin{align*}
+\left(\alpha-3\gamma+\gamma\left[\left(1-\frac{1}{\delta}\right)\left(\frac{zF_\delta'(z)}{F_\delta(z)}\right)+
\frac{1}{\delta}\left(1+\frac{zF_\delta''(z)}{F_\delta'(z)}\right)\right]\right)
\left(\frac{F_\delta(z)}{z}\right)^\delta\left(\frac{zF_\delta'(z)}{F_\delta(z)}\right)
\end{align*}}}
Therefore, from \eqref{eq-F/z_Series}, a simple calculation
gives
\begin{align*}
H_0(z)=1+2(1-\beta)\sum_{n=1}^\infty \tau_nz^n.
\end{align*}
Now rewriting \eqref{eq-generalized-betacond-main2} as
\begin{align*}
\dfrac{1}{1-\beta}=-\dfrac{2}{(1-\xi)}\sum_{n=1}^\infty (-1)^n\tau_n.
\end{align*}
and finding $H_0(z)$ at $z=-1$  gives $H_0(-1)= =\xi$, which clearly implies the sharpness
of the result.
\qed

\noindent
{\underline{\bf{Proof of Theorem $\ref{thm-univ:hyper}$}}}.
Since $\mathcal{H}^\delta=z^\delta{\,}_2F_1(a,b;c;z)\ast
f(z)$. it is a simple exercise to note that
\begin{align*}
\left(\dfrac{\mathcal{H}}{z}\right)^\delta\left(\dfrac{z\mathcal{H}'}{\mathcal{H}}\right)
=N_1(z)\ast\left(\dfrac{f(z)}{z}\right)^\delta,
\end{align*}
where
$N_1(z):=\,_2F_1(a,b;c;z)+\frac{ab}{c\delta}z\,\,_2F_1(a+1,b+1;c+1;z)$.
Further from the above expression, it follows that
\begin{align*}
\left[\left(1-\dfrac{1}{\delta}\right)\left(\dfrac{z\mathcal{H}'}{\mathcal{H}}\right)+\dfrac{1}{\delta}\left(\dfrac{z\mathcal{H}''}{\mathcal{H'}}\right)\right]
\left(\dfrac{\mathcal{H}}{z}\right)^\delta\left(\dfrac{z\mathcal{H}'}{\mathcal{H}}\right)
=N_2(z)\ast\left(\dfrac{f(z)}{z}\right)^\delta,
\end{align*}
where
\begin{align*}
N_2(z):=\left(1-\dfrac{1}{\delta}\right)\,_2F_1(a,b;c;z)+\frac{2ab}{c\delta}z\,\,_2F_1(a+1,b+1;c+1;z)\\
+\frac{a(a+1)b(b+1)}{c(c+1)\delta^2}z^2\,\,_2F_1(a+2,b+2;c+2;z).
\end{align*}
Therefore
\newline
$\displaystyle
\left(\alpha-3\gamma+\gamma\left[\left(1-\frac{1}{\delta}\right)\left(\frac{zf'}{f}\right)+
\frac{1}{\delta}\left(1+\frac{zf''}{f'}\right)\right]\right)
\left(\frac{f}{z}\right)^\delta\!\left(\frac{zf'}{f}\right)
$
{\small{
\begin{align*}
+(1-\alpha+2\gamma)\left(\frac{f}{z}\right)^\delta
=N_3(z)\ast\left(\dfrac{f(z)}{z}\right)^\delta,
\end{align*}}}
where
\begin{align*}
N_3(z):=\,_2F_1(a,b;c;z)+\frac{ab}{c\delta}\left(\alpha-\gamma+\dfrac{\gamma}{\delta}\right)z\,\,_2F_1(a+1,b+1;c+1;z)\\
+\frac{a(a+1)b(b+1)\gamma}{c(c+1)\delta^2}z^2\,\,_2F_1(a+2,b+2;c+2;z).
\end{align*}
From the contiguous relation for Gaussian hypergeometric
functions \cite[Page 96]{AndrewSpecialFn} (or the verification
can be made by comparing the coefficient of $z^n$ both the
sides)
\begin{align*}
b\,_2F_1(a+1,b+1;c+1;z)=c\left(\,_2F_1(a+1,b;c;z)-\,_2F_1(a,b;c;z)\right),
\end{align*}
it follows that
{\small{
\begin{align*}
N_3(z)=\left[1\!-\!\dfrac{a}{\delta}\left(\alpha\!-\!\gamma\left(1\!+\!\dfrac{a}{\delta}\right)\right)\right]\,_2F_1(a,b;c;z)+
\dfrac{a}{\delta}\left[\alpha\!-\!\gamma\left(1\!+\!\dfrac{(2a+1)}{\delta}\right)\right]\,_2F_1(a+1,b;c;z)\\
+\dfrac{a(a+1)\gamma}{\delta^2}\,_2F_1(a+2,b;c;z).
\end{align*}}}
For $a>0$, $b>0$ and $c>0$, the integral representation of
Gaussian hypergeometric function is given as \cite{KimRonning}

{\footnotesize{
\begin{align*}
_2F_1(a,b;c;z)=
\dfrac{\Gamma(c)}{\Gamma(a)\Gamma(b)\Gamma(c\!-\!a\!-\!b\!+\!1)}\int_0^1\!\!
t^{b-1}(1-t)^{c-a-b}{\,}_{2}F_1\left(\!\!\!\!
\begin{array}{cll}&\displaystyle c-a,\quad 1-a
\\
&\displaystyle c-a-b+1
\end{array};1-t\right)\dfrac{1}{(1-tz)}dt.
\end{align*}}}
In view of the above integral form, it is easy to see that
\begin{align*}
N_3(z)=\dfrac{\Gamma(c)}{\Gamma(a)\Gamma(b)\Gamma(c-a-b-1)}\int_0^1 t^{b-1}(1-t)^{c-a-b-2}\dfrac{N_4(t)}{(1-tz)}dt,
\end{align*}
where
\begin{align*}
N_4(t):=&\left[1-\dfrac{a}{\delta}\left(\alpha-\gamma\left(1+\dfrac{a}{\delta}\right)\right)\right](1-t)^2\,\,\dfrac{_2F_1(c-a,1-a;c-a-b+1;1-t)}
{(c-a-b)(c-a-b-1)}\\
&+\dfrac{a}{\delta}\left[\alpha-\gamma\left(1+\dfrac{(2a+1)}{\delta}\right)\right](1-t)\,\,\dfrac{_2F_1(c-a-1,-a;c-a-b;1-t)}
{a(c-a-b-1)}\\
&+\dfrac{\gamma}{\delta^2}\,_2F_1(c-a-2,-(a+1);c-a-b-1;1-t)\\
=&e_1+e_2(1-t)+\sum_{n=0}^\infty\,e_3\dfrac{(c-a)_n(1-a)_n}{(c-a-b)_{n+1}(3)_n}(1-t)^{n+2},
\end{align*}
where
\begin{align*}
e_1&:=\dfrac{\gamma}{\delta^2}\\
e_2&:=\dfrac{1}{\delta^2(c-a-b-1)}\left(\alpha\delta-\gamma(\delta+2a+1)-\gamma(a+1)(c-a-2)\right)\\
e_3&:=\dfrac{1}{2\delta^2(c-a-b)(c-a-b-1)}\left(n^2(\delta^2-a\alpha\delta+a\gamma\delta+a^2\gamma)\dfrac{}{}\right.\\
&\quad+n\left[3(\delta^2-a\alpha\delta+a\gamma\delta+a^2\gamma)-a(\alpha\delta-\gamma(\delta+2a+1))(c-a-1)\dfrac{}{}\right]\\
&\quad+2(\delta^2-a\alpha\delta+a\gamma\delta+a^2\gamma)-a(c-a-1)\left[\dfrac{}{}2(\alpha\delta-\gamma(\delta+2a+1))\right.\\
&\quad\left.\left.\dfrac{}{}-\gamma(a+1)(c-a-2)\right]\right),
\end{align*}
which is non-negative when $0<1+b<(c-a)<2$,
$0<a\leq\min\left\{1,\frac{\delta}{2(\alpha-\gamma)}\right\}$
and $\alpha>\gamma\left(1+\frac{(2a+1)}{\delta}\right)\geq0$. Since
${\rm Re} \left(\frac{1}{1-tz}\right)>\frac{1}{(1+t)}$, for
$|z|<1$, which on using $N_3(z)$ gives ${\rm
Re}\,N_3(z)>N_3(-1)$. Applying Lemma \ref{Lemma_Univalence}
will give the required result.

Now, to verify the sharpness, consider the function
\begin{align*}
\left(\dfrac{f(z)}{z}\right)^\delta=\beta_1+(1-\beta_1)\dfrac{1+z}{1-z}\quad{\rm and}
\quad N_3(z)=\beta_2+(1-\beta_2)\dfrac{1+z}{1-z}.
\end{align*}
Since $\beta=1-2(1-\beta_1)(1-\beta_2)$, therefore
\begin{align*}
\left(\dfrac{f(z)}{z}\right)^\delta\ast N_3(z)&=1+4(1-\beta_1)(1-\beta_2)\dfrac{1+z}{1-z}
=\beta+(1-\beta)\dfrac{1+z}{1-z}
\end{align*}
which clearly implies the sharpness of the result.
\qed

\end{document}